\documentclass[12pt,a4paper]{amsart}

\usepackage[english]{babel}
\usepackage[T1]{fontenc}
\usepackage[utf8]{inputenc}
\usepackage{amsmath,amsthm,amsfonts,amssymb,mathtools}
\usepackage[colorlinks=true,linkcolor=blue,citecolor=blue,urlcolor=blue]{hyperref}
\usepackage[capitalise,nameinlink]{cleveref}
\usepackage[margin=1.15in]{geometry}

\newcommand{\U}{\mathcal{U}}
\newcommand{\C}{\mathsf{C}}
\newcommand{\Sub}{\mathcal{S}}
\newcommand{\f}{\mathtt{f}}
\newcommand{\HJ}{\operatorname{HJ}}
\newcommand{\HJeq}{\operatorname{HJ}^{=}}
\newcommand{\WC}{\operatorname{WC}}
\newcommand{\Sym}{\operatorname{Sym}}
\newcommand{\cont}{\operatorname{cont}}
\newcommand{\restr}{\upharpoonright}
\newcommand{\N}{\mathbb{N}}
\newcommand{\F}{\mathbb{F}}

\theoremstyle{plain}
\newtheorem{theorem}{Theorem}[section]
\newtheorem{lemma}[theorem]{Lemma}
\newtheorem{corollary}[theorem]{Corollary}
\newtheorem{proposition}[theorem]{Proposition}

\theoremstyle{definition}
\newtheorem{definition}[theorem]{Definition}
\newtheorem{counterexample}[theorem]{Counterexample}

\theoremstyle{remark}
\newtheorem{remark}[theorem]{Remark}

\title[From Shelah's block-content to Hales--Jewett]{From Shelah's block-content to Hales--Jewett}

\author{Mohammad Golshani}
\address{School of Mathematics, Institute for Research in Fundamental Sciences (IPM), P.O. Box 19395--5746, Tehran, Iran}
\email{golshani.m@gmail.com}

\author{Mostafa Mirabi}
\address{The Taft School, Watertown, CT 06795, USA, and Wesleyan University, Middletown, CT 06459, USA}
\email{mmirabi@wesleyan.edu}

\author{Saharon Shelah}
\address{Einstein Institute of Mathematics, The Hebrew University of Jerusalem, Jerusalem 91904, Israel, and Department of Mathematics, Rutgers University, New Brunswick, NJ 08854, USA}
\email{shelah@math.huji.ac.il}

\thanks{The first author's research has been supported by a grant from IPM (No. 1405030417). The third author research partially supported by the Israel Science Foundation
(ISF) grant no: 1838/19, and Israel Science Foundation (ISF) grant no: 2320/23; Research partially
supported by the grant “Independent Theories” NSF-BSF, (BSF 3013005232). This is publication
number 1182 of the third author.}


\begin{document}
\begin{abstract}
We study the quantitative relationship between the Hales--Jewett numbers and
Shelah's block-content canonization functions. Block-content canonization yields
a block subspace on which the color of a word is determined solely by the
multiplicities of the alphabet letters among the variable blocks. We show that
this canonical information, combined with the multidimensional Gallai--Witt
theorem, suffices to produce a monochromatic Hales--Jewett subspace. The
argument passes to the space of content vectors, finds a monochromatic
homothetic copy of a finite content simplex, and lifts it through the canonical
block subspace. Combined with the elementary fact that Hales--Jewett bounds
block canonization, this gives a two-way quantitative comparison up to an
explicit change of parameters. The underlying mechanism may be summarized by
the slogan
\begin{center}
Hales--Jewett = block canonization + Gallai--Witt.
\end{center}
We also prove the corresponding equal-block result and show, by an explicit
coloring over every finite field of odd prime order, that the analogous
singleton-coordinate canonization principle fails as soon as two coordinates
remain live.
\end{abstract}

\maketitle

\section{Introduction}

The Hales--Jewett theorem is a central result of finite Ramsey theory. It asserts that, for every finite alphabet and every finite coloring, a sufficiently high-dimensional combinatorial cube contains a monochromatic combinatorial subspace \cite{hj}. For a finite alphabet $\Lambda$ and a finite color set $\C$, we write $\HJ_{\C}(n,\Lambda)$ for the least dimension that guarantees a monochromatic $n$-dimensional subspace.

Shelah introduced a family of finite partition functions designed to isolate intermediate Ramsey and canonization phenomena \cite{shelah-model,shelah-partition}. In this paper we study the block-content functions $\f^{8,*}$ and $\f^{9,*}$. Informally, $\f^{8,*}_{\Lambda}(m,\C)$ is the least dimension in which every $\C$-coloring admits a block subspace with $m$ variable blocks such that the induced color is invariant under permutations of those blocks. Equivalently, the color depends only on the content vector
\[
\cont_{\Lambda}(\eta)
=
\bigl\langle |\eta^{-1}\{\alpha\}|:\alpha\in\Lambda\bigr\rangle.
\]
This conclusion is weaker than monochromaticity: different content vectors may still receive different colors. Our main result shows that a homothetic Ramsey argument in the space of content vectors removes this remaining dependence.

Let
\[
\Delta_h(M)=\left\{\vec{x}\in\N^h:\sum_{e<h}x_e\leq M\right\}.
\]
The points of $\Delta_h(M)$ may be viewed as the multiplicity vectors of words of length $M$ over an alphabet of size $h+1$, with the last multiplicity determined by the first $h$. Let $S_{\C}(h,n)$ denote the least $M$ such that every coloring of $\Delta_h(M)$ contains a monochromatic homothetic copy of $\Delta_h(n)$. The main comparison is
\begin{equation}\label{eq:intro-main}
\HJ_{\C}(n,\Lambda)
\leq
\f^{8,*}_{\Lambda}\bigl(S_{\C}(|\Lambda|-1,n),\C\bigr),
\qquad |\Lambda|>1.
\end{equation}
The multiset parameter is controlled by the multidimensional Gallai--Witt number:
\begin{equation}\label{eq:intro-gw}
S_{\C}(h,n)
\leq
h\bigl(\WC_{\C}(h,n)-1\bigr).
\end{equation}
Combining \eqref{eq:intro-main} and \eqref{eq:intro-gw} gives
\begin{equation}\label{eq:intro-consequence}
\HJ_{\C}(n,\Lambda)
\leq
\f^{8,*}_{\Lambda}\Bigl((|\Lambda|-1)\bigl(\WC_{\C}(|\Lambda|-1,n)-1\bigr),\C\Bigr).
\end{equation}

The proof has three steps. First, $\f^{8,*}$ passes from an arbitrary coloring to one that depends only on block content. Second, the resulting coloring is regarded as a coloring of the simplex $\Delta_{|\Lambda|-1}(M)$. Third, a monochromatic homothetic copy of $\Delta_{|\Lambda|-1}(n)$ is lifted through the original block embedding to a monochromatic $n$-dimensional combinatorial subspace. The appearance of $|\Lambda|-1$ rather than $|\Lambda|$ is intrinsic: once the total number of variable blocks is fixed, the multiplicity of the final letter is determined by the others.

The comparison also has a reverse direction at the same parameter. A monochromatic $m$-dimensional subspace is automatically block-content canonical, and hence
\[
\f^{8,*}_{\Lambda}(m,\C)
\leq
\HJ_{\C}(m,\Lambda).
\]
Thus $\f^{8,*}$ and the Hales--Jewett numbers bound one another: the same-parameter inequality goes from Hales--Jewett to canonization, while the converse uses the content-simplex parameter, which is explicitly controlled by Gallai--Witt through \eqref{eq:intro-consequence}. We establish the parallel result for equal-size variable blocks:
\[
\HJeq_{\C}(n,\Lambda)
\leq
\f^{9,*}_{\Lambda}\bigl(S_{\C}(|\Lambda|-1,n),\C\bigr).
\]

The use of blocks is essential. We introduce the natural singleton-coordinate analogue of block-content canonization and exhibit, for every odd prime $p$, a coloring of $\F_p^M$ that has no singleton-content-canonical set of size two. In particular, the block structure in $\f^{8,*}$ is not merely a convenient reformulation; it is precisely what makes the comparison with Hales--Jewett possible.

The paper is organized as follows. In \Cref{sec:preliminaries} we fix notation, define the relevant Ramsey numbers, and collect elementary comparisons. In \Cref{sec:content} we reformulate block canonization in terms of content and prove the singleton-coordinate obstruction. The main estimates are proved in \Cref{sec:main}, and their consequences are summarized in \Cref{sec:comparisons}.

\section{Preliminaries}\label{sec:preliminaries}

Throughout, $i,j,k,\ell,m,n$ denote natural numbers, and \emph{positive} means nonzero. We identify each natural number $n$ with the ordered set $\{0,\ldots,n-1\}$. The alphabet $\Lambda$ and the color set $\C$ are always finite and nonempty.

All Ramsey numbers below are understood as least integers satisfying the stated property, when such integers exist. An inequality between two such numbers means that whenever the right-hand side exists, so does the left-hand side, and the left-hand side is no larger. Equivalently, the inequalities may be read in the extended natural numbers, with a nonexistent value interpreted as $+\infty$.

\subsection{Combinatorial subspaces and Hales--Jewett numbers}

If $M$ is a finite linear order, let
\[
\U_{M,\Lambda}=\{\eta:M\to\Lambda\}.
\]
When $M=n$ with its usual order, we write $\U_{n,\Lambda}$.

\begin{definition}\label{def:block-subspace}
Let $B_0,\ldots,B_{m-1}\subseteq M$ be pairwise disjoint nonempty sets, and let
\[
\rho:M\setminus\bigcup_{i<m}B_i\longrightarrow\Lambda.
\]
The associated $m$-dimensional block subspace is
\[
\Sub(\langle B_i:i<m\rangle,\rho)
=
\left\{\nu\in\U_{M,\Lambda}:
\nu\supseteq\rho\text{ and }\nu\restr B_i\text{ is constant for every }i<m
\right\}.
\]
Its block embedding is the map
\[
\Phi_{\langle B_i:i<m\rangle,\rho}:\U_{m,\Lambda}\longrightarrow\U_{M,\Lambda}
\]
defined by
\[
\Phi_{\langle B_i:i<m\rangle,\rho}(\eta)(a)=
\begin{cases}
\eta(i),&a\in B_i,\\
\rho(a),&a\notin\displaystyle\bigcup_{i<m}B_i.
\end{cases}
\]
The block subspace is exactly the image of this embedding. When $m=1$, it is called a combinatorial line.
\end{definition}

For $\eta\in\U_{m,\Lambda}$, define its \emph{content} by
\[
\cont_{\Lambda}(\eta)
=
\bigl\langle |\eta^{-1}\{\alpha\}|:\alpha\in\Lambda\bigr\rangle.
\]

\begin{definition}\label{def:HJ}
For a positive integer $m$, $\HJ_{\C}(m,\Lambda)$ is the least $k$, if it exists, such that for every finite linear order $M$ of size $k$ and every coloring
\[
\mathbf d:\U_{M,\Lambda}\longrightarrow\C,
\]
there is a $\mathbf d$-monochromatic $m$-dimensional block subspace of $\U_{M,\Lambda}$.
\end{definition}

\begin{definition}\label{def:HJeq}
For a positive integer $m$, $\HJeq_{\C}(m,\Lambda)$ is the least $k$, if it exists, such that every coloring of $\U_{M,\Lambda}$, with $|M|=k$, has a monochromatic $m$-dimensional block subspace
\[
\Sub(\langle B_i:i<m\rangle,\rho)
\]
whose variable blocks all have the same size.
\end{definition}

The Hales--Jewett theorem \cite{hj} states that $\HJ_{\C}(m,\Lambda)$ is finite. Shelah's proof gives primitive recursive upper bounds \cite{shelah-vdw}; see also the expositions in \cite{nilli,matet}.

\subsection{Gallai--Witt numbers and content simplices}

\begin{definition}\label{def:GW}
Let $h$ and $m$ be positive integers. The number $\WC_{\C}(h,m)$ is the least $N$, if it exists, such that every coloring
\[
\mathbf d:\U_{h,N}\longrightarrow\C
\]
admits $d>0$ and $a_0,\ldots,a_{h-1}<N$ satisfying
\[
a_e+dm<N\qquad(e<h)
\]
and such that $\mathbf d$ is constant on the grid
\[
\left\{\langle a_e+di_e:e<h\rangle:i_0,\ldots,i_{h-1}\leq m\right\}.
\]
\end{definition}

The finiteness of $\WC_{\C}(h,m)$ is the multidimensional Gallai--Witt theorem \cite{rado,witt,graham}. For $h=1$, it is van der Waerden's theorem \cite{waerden}.

For $h,M\in\N$, set
\[
\Delta_h(M)=\left\{\vec{x}\in\N^h:\sum_{e<h}x_e\leq M\right\}.
\]
We regard $\Delta_h(M)$ as the set of multiplicity vectors of words of length $M$ over $h+1$ symbols: the first $h$ coordinates record the first $h$ multiplicities, while the final multiplicity is determined by their sum.

\begin{definition}\label{def:S}
Let $S_{\C}(h,n)$ be the least $M$, if it exists, such that every coloring
\[
c:\Delta_h(M)\longrightarrow\C
\]
contains a monochromatic homothetic copy of $\Delta_h(n)$, that is, a set
\[
\vec{a}+d\Delta_h(n)
=
\{\vec{a}+d\vec{y}:\vec{y}\in\Delta_h(n)\}
\]
for some $d>0$ and some $\vec{a}\in\N^h$.
\end{definition}

\begin{proposition}\label{prop:S-WC}
For all positive integers $h,n$,
\[
S_{\C}(h,n)
\leq
h\bigl(\WC_{\C}(h,n)-1\bigr).
\]
\end{proposition}

\begin{proof}
Let $N=\WC_{\C}(h,n)$ and put $M=h(N-1)$. Given a coloring
\[
c:\Delta_h(M)\longrightarrow\C,
\]
restrict it to the cube $\{0,\ldots,N-1\}^h$, which is contained in $\Delta_h(M)$. By the definition of $\WC_{\C}(h,n)$, there are $d>0$ and $a_0,\ldots,a_{h-1}<N$ such that $a_e+dn<N$ for every $e<h$ and $c$ is constant on
\[
\left\{\langle a_e+di_e:e<h\rangle:i_0,\ldots,i_{h-1}\leq n\right\}.
\]
Since $\Delta_h(n)\subseteq\{0,\ldots,n\}^h$, the homothetic copy
\[
\vec{a}+d\Delta_h(n)
\]
is contained in this monochromatic grid. Hence it is monochromatic.
\end{proof}

\subsection{Shelah's block partition functions}

We now recall the four partition functions used in the paper. The balanced functions $\f^8$ and $\f^9$ require $|\Lambda|$ to divide the number of variable blocks; the canonization functions $\f^{8,*}$ and $\f^{9,*}$ are defined for every positive number of variable blocks.

\begin{definition}\label{def:f8f9}
Assume that $m$ is positive and that $|\Lambda|$ divides $m$. The number $\f^9_{\Lambda}(m,\C)$ is the least $k$, if it exists, such that for every finite linear order $M$ of size $k$ and every coloring $\mathbf d:\U_{M,\Lambda}\to\C$, there are pairwise disjoint nonempty blocks $B_i\subseteq M$, $i<m$, and a map
\[
\rho:M\setminus\bigcup_{i<m}B_i\longrightarrow\Lambda
\]
such that:
\begin{enumerate}
\item $|B_i|=|B_j|$ for all $i,j<m$;
\item $\mathbf d$ is constant on the set of all $\nu\in\Sub(\langle B_i:i<m\rangle,\rho)$ satisfying
\[
\bigl|\{i<m:\nu\restr B_i\equiv\alpha\}\bigr|
=
\frac{m}{|\Lambda|}
\]
for every $\alpha\in\Lambda$.
\end{enumerate}
The function $\f^8_{\Lambda}(m,\C)$ is defined in the same way, but without the equal-size requirement in item~(1).
\end{definition}

\begin{definition}\label{def:f8starf9star}
Let $m$ be positive. The number $\f^{9,*}_{\Lambda}(m,\C)$ is the least $k$, if it exists, such that for every finite linear order $M$ of size $k$ and every coloring $\mathbf d:\U_{M,\Lambda}\to\C$, there are pairwise disjoint nonempty blocks $B_i\subseteq M$, $i<m$, and a map
\[
\rho:M\setminus\bigcup_{i<m}B_i\longrightarrow\Lambda
\]
such that:
\begin{enumerate}
\item $|B_i|=|B_j|$ for all $i,j<m$;
\item whenever $\nu_1,\nu_2\in\Sub(\langle B_i:i<m\rangle,\rho)$ satisfy
\[
\bigl|\{i<m:\nu_1\restr B_i\equiv\alpha\}\bigr|
=
\bigl|\{i<m:\nu_2\restr B_i\equiv\alpha\}\bigr|
\]
for every $\alpha\in\Lambda$, then $\mathbf d(\nu_1)=\mathbf d(\nu_2)$.
\end{enumerate}
The function $\f^{8,*}_{\Lambda}(m,\C)$ is defined in the same way, but without the equal-size requirement in item~(1).
\end{definition}

\begin{lemma}\label{lem:monotonicity}
For fixed finite $\Lambda$ and $\C$, each of the functions
\begin{align*}
m&\longmapsto\HJ_{\C}(m,\Lambda),
& m&\longmapsto\HJeq_{\C}(m,\Lambda),\\
m&\longmapsto\f^{8,*}_{\Lambda}(m,\C),
& m&\longmapsto\f^{9,*}_{\Lambda}(m,\C)
\end{align*}
is nondecreasing.
\end{lemma}

\begin{proof}
Let $m'\leq m$. From an $m$-dimensional witness, retain any $m'$ variable blocks and freeze the remaining blocks. This gives the required $m'$-dimensional witness. For the canonization functions, two assignments with the same content on the retained blocks have the same full content after the discarded blocks are frozen. The equal-block property is preserved when blocks are discarded.
\end{proof}

\begin{lemma}\label{lem:basic-comparisons}
The following pointwise inequalities hold whenever the quantities involved are defined:
\begin{enumerate}
\item $\f^8_{\Lambda}(m,\C)\leq\f^9_{\Lambda}(m,\C)$;
\item $\f^{8,*}_{\Lambda}(m,\C)\leq\f^{9,*}_{\Lambda}(m,\C)$;
\item $\f^\ell_{\Lambda}(m,\C)\leq\f^{\ell,*}_{\Lambda}(m,\C)$ for $\ell\in\{8,9\}$, whenever $|\Lambda|$ divides $m$;
\item $\f^{8,*}_{\Lambda}(m,\C)\leq\HJ_{\C}(m,\Lambda)$;
\item $\f^{9,*}_{\Lambda}(m,\C)\leq\HJeq_{\C}(m,\Lambda)$.
\end{enumerate}
\end{lemma}

\begin{proof}
The first two inequalities follow by dropping the equal-size requirement on the variable blocks. For item~(3), content canonization implies monochromaticity on the balanced part of the subspace, because all balanced words have the same content vector. Items~(4) and~(5) follow because a monochromatic subspace, respectively a monochromatic equal-block subspace, satisfies the corresponding canonization condition.
\end{proof}

\begin{lemma}\label{lem:HJeq-bound}
For every positive integer $m$,
\[
\HJeq_{\C}(m,\Lambda)
\leq
m\cdot\HJ_{\C}(1,\Lambda^m).
\]
Consequently,
\[
\f^{9,*}_{\Lambda}(m,\C)
\leq
m\cdot\HJ_{\C}(1,\Lambda^m),
\]
where $\Lambda^m$ is regarded as a finite alphabet.
\end{lemma}

\begin{proof}
Let $\Gamma=\Lambda^m$ and $N=\HJ_{\C}(1,\Gamma)$. Up to order isomorphism, it is enough to work with the lexicographically ordered set $M=N\times m$. Given a coloring $\mathbf d:\U_{M,\Lambda}\to\C$, define
\[
F:\U_{N,\Gamma}\longrightarrow\U_{M,\Lambda}
\]
by
\[
F(\eta)(i,\ell)=\eta(i)(\ell).
\]
Set $\mathbf e(\eta)=\mathbf d(F(\eta))$. By the definition of $N$, there are a nonempty set $I\subseteq N$ and a map $\rho:N\setminus I\to\Gamma$ such that $\mathbf e$ is constant on the corresponding combinatorial line in $\U_{N,\Gamma}$.

For each $\ell<m$, put
\[
B_\ell=I\times\{\ell\}.
\]
Define $\varrho:M\setminus\bigcup_{\ell<m}B_\ell\to\Lambda$ by
\[
\varrho(i,\ell)=\rho(i)(\ell)
\qquad(i\in N\setminus I,\ \ell<m).
\]
The blocks $B_0,\ldots,B_{m-1}$ are nonempty and have equal size, and
\[
\Sub(\langle B_\ell:\ell<m\rangle,\varrho)
\]
is the image under $F$ of the monochromatic line. It is therefore $\mathbf d$-monochromatic. The second inequality follows from Lemma~\ref{lem:basic-comparisons}.
\end{proof}

\section{Canonization}\label{sec:content}

We first restate the defining property of $\f^{8,*}$ in a form suited to the proof of the main theorem.

Let $E_m$ be the equivalence relation on $\U_{m,\Lambda}$ induced by the full symmetric group on $m$:
\[
\eta_1E_m\eta_2
\quad\Longleftrightarrow\quad
\exists\pi\in\Sym(m)\ \bigl(\eta_2=\eta_1\circ\pi\bigr).
\]
Equivalently,
\[
\eta_1E_m\eta_2
\quad\Longleftrightarrow\quad
\cont_{\Lambda}(\eta_1)=\cont_{\Lambda}(\eta_2).
\]

\begin{proposition}\label{prop:content-form-f8star}
Let $m$ be positive and $k=\f^{8,*}_{\Lambda}(m,\C)$. For every finite linear order $M$ of size $k$ and every coloring $\mathbf d:\U_{M,\Lambda}\to\C$, there are pairwise disjoint nonempty blocks $B_i\subseteq M$, $i<m$, and a map
\[
\rho:M\setminus\bigcup_{i<m}B_i\longrightarrow\Lambda
\]
such that the induced coloring
\[
D:\U_{m,\Lambda}\longrightarrow\C,
\qquad
D(\eta)=\mathbf d\bigl(\Phi_{\langle B_i:i<m\rangle,\rho}(\eta)\bigr),
\]
is invariant under $E_m$. Equivalently, $D$ depends only on the content of $\eta$.
\end{proposition}

\begin{proof}
This is the defining property of $\f^{8,*}_{\Lambda}(m,\C)$ written on the standard cube $\U_{m,\Lambda}$. Two words in $\U_{m,\Lambda}$ have the same content exactly when the corresponding points of the block subspace assign each letter to the same number of variable blocks.
\end{proof}

\subsection{Failure of singleton-coordinate canonization}\label{sec:singletons}

The block formulation in $\f^{8,*}$ is essential rather than notational. Consider the following natural singleton-coordinate analogue.

\begin{definition}\label{def:f13sing}
Let $\f^{13,\mathrm{sing}}_{\Lambda}(m,\C)$ be the least $k$, if it exists, such that for every finite linear order $M$ of size $k$ and every coloring $\mathbf d:\U_{M,\Lambda}\to\C$, there are a set $A\subseteq M$ of size $m$ and a map $\rho:M\setminus A\to\Lambda$ such that, for all $\eta_1,\eta_2\in\U_{A,\Lambda}$,
\[
\cont_{\Lambda}(\eta_1)=\cont_{\Lambda}(\eta_2)
\quad\Longrightarrow\quad
\mathbf d(\eta_1\cup\rho)=\mathbf d(\eta_2\cup\rho).
\]
\end{definition}

Thus $\f^{13,\mathrm{sing}}$ asks for content canonization on a set of singleton coordinates rather than on a family of variable blocks.

\begin{counterexample}\label{counter:singletons}
Let $p$ be an odd prime. Then, for every $m\geq2$,
\[
\f^{13,\mathrm{sing}}_{\F_p}(m,\F_p)=\infty.
\]
\end{counterexample}

\begin{proof}
For every finite linear order $M$, define a coloring
\[
d:\U_{M,\F_p}\longrightarrow\F_p
\]
by
\[
d(\eta)=\sum_{a<_M b}\eta(a)\eta(b)^2.
\]
We show that no set of at least two live coordinates can be singleton-content canonical.

Suppose that a set $A$ of at least two live coordinates and an outside word $\rho:M\setminus A\to\F_p$ witness singleton-content canonization. Choose $u<v$ in $A$ and set every coordinate in $A\setminus\{u,v\}$ equal to $0$. The resulting two-variable coloring has the form
\[
D(x,y)=xy^2+P(x)+Q(y)+c,
\]
where $P,Q:\F_p\to\F_p$ and $c\in\F_p$. Since the assignments $(x,y)$ and $(y,x)$ have the same content, canonization gives
\[
D(x,y)=D(y,x)
\]
for all $x,y\in\F_p$. Hence
\[
xy^2-yx^2=P(y)-P(x)+Q(x)-Q(y).
\]
Setting $x=0$ yields
\[
Q(y)-Q(0)=P(y)-P(0)
\]
for every $y\in\F_p$, so $P-Q$ is constant. The right-hand side of the preceding identity is therefore always $0$, and consequently
\[
xy^2-yx^2=0
\]
for all $x,y\in\F_p$. Taking $x=1$ and $y=2$ gives $2=0$ in $\F_p$, a contradiction because $p$ is odd. Thus no singleton-content-canonical live set of size at least two exists.
\end{proof}

\begin{remark}
The obstruction explains why the block version is the appropriate object for comparison with Hales--Jewett. Hales--Jewett variables are represented by blocks of coordinates, exactly the setting of $\f^{8,*}$ used in \Cref{thm:main-HJ-f8star}.
\end{remark}

\section{The main bounds}\label{sec:main}

We now combine block-content canonization with the content-simplex form of the Gallai--Witt theorem.

\begin{theorem}\label{thm:main-HJ-f8star}
Let $\Lambda$ be a finite alphabet with $|\Lambda|>1$, let $\C$ be a finite color set, and let $n$ be positive. Then
\[
\HJ_{\C}(n,\Lambda)
\leq
\f^{8,*}_{\Lambda}\bigl(S_{\C}(|\Lambda|-1,n),\C\bigr).
\]
Consequently,
\[
\HJ_{\C}(n,\Lambda)
\leq
\f^{8,*}_{\Lambda}\Bigl((|\Lambda|-1)\bigl(\WC_{\C}(|\Lambda|-1,n)-1\bigr),\C\Bigr).
\]
\end{theorem}

\begin{proof}
Let $r=|\Lambda|$ and enumerate
\[
\Lambda=\{\alpha_0,\ldots,\alpha_{r-1}\}.
\]
Put
\[
M_0=S_{\C}(r-1,n)
\qquad\text{and}\qquad
K=\f^{8,*}_{\Lambda}(M_0,\C).
\]
Let $M$ be a finite linear order of size $K$, and let
\[
\mathbf c:\U_{M,\Lambda}\longrightarrow\C
\]
be a coloring. By Proposition~\ref{prop:content-form-f8star}, there are pairwise disjoint nonempty blocks $B_j\subseteq M$, $j<M_0$, and a map
\[
\rho^*:M\setminus\bigcup_{j<M_0}B_j\longrightarrow\Lambda
\]
such that the induced coloring
\[
D:\U_{M_0,\Lambda}\longrightarrow\C,
\qquad
D(\eta)=\mathbf c\bigl(\Phi_{\langle B_j:j<M_0\rangle,\rho^*}(\eta)\bigr),
\]
depends only on the content of $\eta$.

The coloring $D$ induces a coloring
\[
D_0:\Delta_{r-1}(M_0)\longrightarrow\C.
\]
For $\vec{x}=(x_0,\ldots,x_{r-2})\in\Delta_{r-1}(M_0)$, choose any $\eta\in\U_{M_0,\Lambda}$ such that
\[
|\eta^{-1}\{\alpha_e\}|=x_e
\qquad(e<r-1),
\]
and define $D_0(\vec{x})=D(\eta)$. This is well-defined: the multiplicity of $\alpha_{r-1}$ is $M_0-\sum_{e<r-1}x_e$, and $D$ depends only on the complete content vector.

By the definition of $S_{\C}(r-1,n)$, there are $d>0$ and
\[
\vec{a}=(a_0,\ldots,a_{r-2})\in\N^{r-1}
\]
such that
\[
\vec{a}+d\Delta_{r-1}(n)
\]
is monochromatic for $D_0$. Since this homothetic copy is contained in $\Delta_{r-1}(M_0)$,
\[
\sum_{e<r-1}a_e+nd\leq M_0.
\]
Choose pairwise disjoint sets of block indices
\[
P_e\subseteq M_0\quad(e<r-1),
\qquad
V_i\subseteq M_0\quad(i<n),
\]
with
\[
|P_e|=a_e
\qquad\text{and}\qquad
|V_i|=d.
\]

For each $v\in\U_{n,\Lambda}$, define $\eta_v\in\U_{M_0,\Lambda}$ by
\[
\eta_v(j)=
\begin{cases}
\alpha_e,&j\in P_e\text{ for some }e<r-1,\\
v(i),&j\in V_i\text{ for some }i<n,\\
\alpha_{r-1},&\text{otherwise}.
\end{cases}
\]
For each $i<n$, set
\[
C_i=\bigcup_{j\in V_i}B_j.
\]
The sets $C_0,\ldots,C_{n-1}$ are pairwise disjoint nonempty blocks. Freeze the original blocks indexed by $P_e$ to $\alpha_e$, freeze all remaining original blocks not used by the $V_i$ to $\alpha_{r-1}$, and retain $\rho^*$ outside the original block subspace. With this frozen word, the map
\[
v\longmapsto
\Phi_{\langle B_j:j<M_0\rangle,\rho^*}(\eta_v)
\]
is the block embedding associated with $C_0,\ldots,C_{n-1}$.

For every $v\in\U_{n,\Lambda}$, the first $r-1$ coordinates of $\cont_{\Lambda}(\eta_v)$ are
\[
\bigl\langle a_e+d|v^{-1}\{\alpha_e\}|:e<r-1\bigr\rangle.
\]
Since
\[
\bigl\langle |v^{-1}\{\alpha_e\}|:e<r-1\bigr\rangle
\in
\Delta_{r-1}(n),
\]
this vector belongs to the monochromatic set $\vec{a}+d\Delta_{r-1}(n)$. Hence all words $\eta_v$ receive the same color under $D$, and their images form a monochromatic $n$-dimensional block subspace of $\U_{M,\Lambda}$.

The second displayed inequality follows from Proposition~\ref{prop:S-WC} and Lemma~\ref{lem:monotonicity}.
\end{proof}

The same argument preserves equal block sizes.

\begin{theorem}\label{thm:main-HJ-f9star}
Let $\Lambda$ be a finite alphabet with $|\Lambda|>1$, let $\C$ be a finite color set, and let $n$ be positive. Then
\[
\HJeq_{\C}(n,\Lambda)
\leq
\f^{9,*}_{\Lambda}\bigl(S_{\C}(|\Lambda|-1,n),\C\bigr).
\]
Consequently,
\[
\HJeq_{\C}(n,\Lambda)
\leq
\f^{9,*}_{\Lambda}\Bigl((|\Lambda|-1)\bigl(\WC_{\C}(|\Lambda|-1,n)-1\bigr),\C\Bigr).
\]
\end{theorem}

\begin{proof}
Repeat the proof of \Cref{thm:main-HJ-f8star} using $\f^{9,*}$ in place of $\f^{8,*}$. The original blocks $B_j$ then have a common size, say $s$. Each final variable block has the form
\[
C_i=\bigcup_{j\in V_i}B_j,
\]
and $|V_i|=d$ for every $i<n$. Thus $|C_i|=ds$ for every $i<n$, so the resulting monochromatic subspace has equal-size variable blocks. The second inequality again follows from Proposition~\ref{prop:S-WC} and Lemma~\ref{lem:monotonicity}.
\end{proof}

\section{Comparisons and consequences}\label{sec:comparisons}

Taking $n=1$ in \Cref{thm:main-HJ-f8star} gives the following line version.

\begin{corollary}\label{cor:line-bound}
For every finite alphabet $\Lambda$ with $|\Lambda|>1$ and every finite color set $\C$,
\[
\HJ_{\C}(1,\Lambda)
\leq
\f^{8,*}_{\Lambda}\Bigl((|\Lambda|-1)\bigl(\WC_{\C}(|\Lambda|-1,1)-1\bigr),\C\Bigr).
\]
\end{corollary}

For fixed finite $\Lambda$ and $\C$, the same-parameter comparisons in Lemma~\ref{lem:basic-comparisons} give
\[
\f^8_{\Lambda}(m,\C)
\leq
\f^{8,*}_{\Lambda}(m,\C)
\leq
\HJ_{\C}(m,\Lambda)
\]
whenever $|\Lambda|$ divides $m$, and
\[
\f^8_{\Lambda}(m,\C)
\leq
\f^9_{\Lambda}(m,\C)
\leq
\f^{9,*}_{\Lambda}(m,\C)
\leq
\HJeq_{\C}(m,\Lambda).
\]
Moreover,
\[
\f^{8,*}_{\Lambda}(m,\C)
\leq
\f^{9,*}_{\Lambda}(m,\C).
\]
Together with \Cref{thm:main-HJ-f8star,thm:main-HJ-f9star}, these inequalities give two-way comparisons between the canonization functions and the corresponding Hales--Jewett numbers. Any upper bound for a Hales--Jewett number immediately bounds the associated canonization number. Conversely, an upper bound for $\f^{8,*}$ or $\f^{9,*}$ yields an upper bound for the corresponding Hales--Jewett number after the content-simplex parameter change, and hence after the explicit Gallai--Witt parameter change in Proposition~\ref{prop:S-WC}.

Known primitive recursive bounds for Hales--Jewett therefore give primitive recursive bounds for $\f^{8,*}$, while Lemma~\ref{lem:HJeq-bound} gives such bounds for $\f^{9,*}$. In the reverse direction, \Cref{thm:main-HJ-f8star} converts any explicit bound for $\f^{8,*}$ into an explicit Hales--Jewett bound after composition with a bound for the multidimensional Gallai--Witt numbers. These are pointwise comparisons of finite numbers; by themselves they do not determine the exact Grzegorczyk complexity of the functions involved.


\begin{thebibliography}{99}

\bibitem{graham}
Ronald L. Graham, Bruce L. Rothschild, and Joel H. Spencer,
\emph{Ramsey Theory}, 2nd ed., John Wiley \& Sons, New York, 1990.

\bibitem{hj}
Alfred W. Hales and Robert I. Jewett,
Regularity and positional games,
\emph{Trans. Amer. Math. Soc.} \textbf{106} (1963), 222--229.

\bibitem{matet}
Pierre Matet,
Shelah's proof of the Hales--Jewett theorem revisited,
\emph{European J. Combin.} \textbf{28} (2007), 1742--1745.

\bibitem{nilli}
Alon Nilli,
Shelah's proof of the Hales--Jewett theorem,
in \emph{Mathematics of Ramsey Theory}, Algorithms and Combinatorics, vol.~5, Springer, Berlin, 1990, pp.~150--151.

\bibitem{rado}
Richard Rado,
Note on combinatorial analysis,
\emph{Proc. London Math. Soc.} \textbf{48} (1943), 122--160.

\bibitem{shelah-vdw}
Saharon Shelah,
Primitive recursive bounds for van der Waerden numbers,
\emph{J. Amer. Math. Soc.} \textbf{1} (1988), 683--697.

\bibitem{shelah-model}
Saharon Shelah,
On what I do not understand (and have something to say), model theory,
\emph{Math. Japon.} \textbf{51} (2000), 329--377.

\bibitem{shelah-partition}
Saharon Shelah,
A partition theorem,
\emph{Sci. Math. Jpn.} \textbf{56} (2002), 413--438.

\bibitem{waerden}
Bartel L. van der Waerden,
Beweis einer Baudetschen Vermutung,
\emph{Nieuw Arch. Wisk.} \textbf{15} (1927), 212--216.

\bibitem{witt}
Ernst Witt,
Ein kombinatorischer Satz der Elementargeometrie,
\emph{Math. Nachr.} \textbf{6} (1951), 261--262.

\end{thebibliography}
\end{document}